\documentclass{article}
\usepackage[utf8]{inputenc}
\usepackage{amssymb}
\usepackage{amsmath}
\newtheorem{lemma}{Lemma}
\newtheorem{theorem}{Theorem}
\newtheorem{corollary}{Corollary}
\title{Average order in wreath products } 
\author{Supravat Sarkar}
\date{}
\begin{document}
\maketitle
\begin{abstract}
    We obtain an exact formula for the average order of elements of
a wreath product of two finite groups. Then focussing our
attention on $p$-groups for primes $p$, we give an estimate for the
average order of a wreath product $A\wr B$ in terms of maximum order of elements of $A$ and
average order of $B$ and an exact formula for the distribution
of orders of elements of $A\wr B.$ Finally, we show how wreath
products can be used to find several rational numbers which are
limits of average orders of a sequence of $p$-groups with
cardinalities going to infinity.
\end{abstract}
\begin{center}
\end{center}
\begin{center}
\end{center}

\section{Introduction}

For a finite group $G,$ define the average order of elements to be
$a(G)=\frac{\sum_{g\in G} order(g)}{|G|}$ and denote the maximum
order of an element of $G$ by $m(G).$ It is well known that among
all groups of a fixed cardinality, the cyclic group has the largest
average order (see \cite{AA, AAI}). Other authors have
studied the average order function with a view to deriving
characterization theorems for nilpotency and solvability (see
\cite{HLM1, HLM2}).

Wreath products of groups often provide examples and
counter-examples for various group theoretic questions. The famous
\textit{Krasner-Kaloujnine embedding theorem} shows that for any two
groups $A$ and $H$, any extension of $A$ by $H$ is isomorphic to a
subgroup of the wreath product $A\wr H$ (see \cite{M}).

In this article, we obtain an exact formula for the average order of
a wreath product of two finite groups. Following that, we focus our
attention on $p$-groups for primes $p$. We give an estimate for the
average order of a wreath product $A\wr B$ in terms of $m(A)$ and
$a(B).$ In addition, we obtain an exact formula for the distribution
of orders of elements of $A\wr B,$ in terms of distributions of
orders of elements of $A$ and $B$. Finally, we show how wreath
products can be used to find several rational numbers which are
limits of average orders of a sequence of $p$-groups with
cardinalities going to infinity.

\section{Notations and Conventions}
\begin{enumerate}
\item Given two groups $A$ and $B$, let $K=\prod_{b\in B} A.$ The group $B$ naturally acts on $K$ by $x\cdot (\alpha_b)_b=(\alpha_{x^{-1}b})_b,$ for $x\in B,$ $(\alpha_b)_b\in K.$ The semidirect product $K\rtimes B$ is called the wreath product of $A$ by $B$ and is denoted by $A\wr B.$
\item $\mu:\mathbf{N}\to \mathbf{Z}$ denotes the Mobius function. Recall that if $n$ is square-free,
$\mu(n)=(-1)^m,$ where $m$ is the number of distinct prime divisors
of $n$. If $n$ is not square-free, $\mu(n)=0.$
\item For a natural number $n,$ define $\tau(n)=\prod_{p|n, p~ prime} (1-p). $ It is easy to see that $\tau(n)=\sum_{d|n} d\mu(d).$
    \item For real valued functions $f$ and $g$ we adopt the notation $g=O(f)$ to mean that there is a constant $M>0$ such that $|g|\leq M|f|$ always.
    \item Define $\psi(A,B)=\frac{a(A\wr B)}{m(A) \cdot a(B)},$ for $p$-groups $A$ and $B.$
    \item Let $A$ be a $p$-group of cardinality $p^a$ and $m(A)=p^d$. For $k\in \mathbf{Z}$ define
    $r_{A,k}=\frac{1}{p^a} \times$ Number of elements of order at most $p^{d-k}.$ \\
    This essentially denotes the cumulative distribution function of the order distribution.  So, $r_{A,k}$ is a non-increasing function of $k,$
    and $r_{A,k}=1$ for $k\leq0,$ $r_{A,k}=0$ for $k>d.$
    \item Fix a prime $p.$ Let us call a real number $\beta$ an \textit{``Average Order Limit"} if there is a sequence of $p$-groups $G_n$ with $|G_n|\to \infty$ and $a(G_n)\to \beta.$

\end{enumerate}
\section{Main Results}
\begin{theorem}
Let $A$, $B$ be finite groups with at least $2$ elements. For
$m\geq1$, let $s_m$ be the number of elements of $A$ whose $m$-th
power is $1$. For $n\vert b$, let $d_n$ be the number of elements of
$B$ of order $n$. Then
$$a(A \wr B)= \sum_{m \bigm| |A|, n \bigm| |B|} \frac{m}{n}
(\frac{s_m}{|A|})^n d_\frac{|B|}{n} \tau(\frac{|A|}{m}).$$
\end{theorem}

Let $p$ be a prime. From now on, assume $A$, $B$ are $p$-groups, $|A|=p^a,$ $|B|=p^b,$ $a,$ $b\geq1.$

\begin{theorem}
 With notations as in theorem 1,
$$a(A \wr B) =p^d a(B) -(p-1)a(B)\sum_{n\leq b} k_n\cdot
[\sum_{m\leq d-1} p^m(\frac{s_{p^m}}{p^a})^{p^n}],$$
 where $k_n=\frac{p^{-n} d_{p^{b-n}}}{a(B)}.$ Note that $k_n\geq 0$ $\forall n$, $\sum_{n\leq b} k_n=1.$
  \end{theorem}

\begin{theorem}
$a(B) \leq a(A \wr B)\leq p^d a(B).$
\end{theorem}
Note that Theorem 3 implies that $0\leq \psi(A,B)\leq 1.$
\begin{theorem} With notations as in theorem 1,
$$a(A \wr B)=p^d a(B) -(p-1)a(B).\sum_{n\leq b}k_n . [ p^{d-1}
(\frac{s_{p^{d-1}}}{p^a})^{p^n}] +O(p^{d-1}a(B)).$$
  \end{theorem}

   Here,  the implicit constant  in $O(p^{d-1}a(B))$ is independent of $p$ as well as of the
   groups. Now, we explicitly write the order distribution of $A\wr B$ in terms of the order distributions of $A$ and $B$.

  \begin{theorem}
    Let m($A$)=$p^d,$  m($B$)=$p^e.$ Then we have m($A\wr B)=p^{d+e},$ and
    $r_{A\wr B,k}=\sum_{i=0}^e (r_{B,i}-r_{B,i+1}) r_{A,k-i}^{p^{b-e+i}}.$
    \end{theorem}
    \begin{corollary}
    $r_{A\wr \mathbf{Z}/p\mathbf{Z},k}=(1-p^{-1})r_{A,k}+\frac{r_{A,k-1}^p}{p}.$
    \end{corollary}

    In view of theorem 4, it is natural to ask whether $a(A\wr B)$
  $=p^d a(B) +O(p^{d-1}a(B))$, i.e., whether ``$\psi(A,B)=1+O(\frac{1}{p})$". We shall show that this is not true in general. However, if we assume $A$ is abelian, then it is true.

    \begin{corollary}
    Let $A$ be a $p$-group. Define $p$-groups $A_n$ recursively by $A_0=A,$ $A_n=A_{n-1}\wr \mathbf{Z}/p\mathbf{Z}.$ Let $B$ be a $p$-group. Then $\psi(A_n,B)\to0$ as $n\to \infty.$
    \end{corollary}

Note that corollary 5.2 shows that ``$\psi(A,B)=1+O(\frac{1}{p})$" cannot be true in general. Now we show how the situation differs if we assume $A$ to be abelian.
\begin{theorem}
For abelian $p$-groups $A,$ define $t(A)$ to be the unique positive
integer such that $A\cong \mathbf{Z}/p^{d_1}\mathbf{Z}\oplus
\mathbf{Z}/p^{d_2}\mathbf{Z}\oplus...\oplus
\mathbf{Z}/p^{d_k}\mathbf{Z}\oplus
(\mathbf{Z}/p^d\mathbf{Z})^{t(A)}$, where $d_1\leq d_2\leq \cdots
\leq d_k<d.$ Then, \\
(i) For each noncyclic group $B,$ $1-p^{-t(A)p}\leq \psi(A,B)\leq
1.$\\
(ii) For cyclic groups $B,$ $\psi(A,B)=1-p^{-t(A)}+O(p^{-t(A)-1}).$\\
In particular, $\psi(A,B)=1+O(\frac{1}{p})$ holds if $A$ is abelian.
\end{theorem}
Now we prove a result about limiting value of average orders for a sequence of groups whose sizes goes to infinity.
\begin{theorem}
Suppose $B_n$ is a sequence of $p$-groups with $a(B_n)\to \beta.$ Then $a((\mathbf{Z}/p\mathbf{Z})^n\wr B_n)\to p\beta.$ So, if $\beta$ is an \textit{Average Order Limit}, then $p\beta$ is also. Hence $p^r\beta$ is an \textit{Average Order Limit} for each nonnegative integer $r$.
\end{theorem}
 Note that whatever the sizes of $B_n$ be, sizes of $(\mathbf{Z}/p\mathbf{Z})^n\wr B_n$ always go to infinity. So,
 taking $\{B_n\}_n$ to be a constant sequence $B$,
we see that $p^r a(B)$ is an \textit{Average Order Limit}, for any
$r\geq2$. This guarantees the existence of many non-integral
\textit{Average Order Limits}. For example, taking
$B=(\mathbf{Z}/p\mathbf{Z})^b$ and $2\leq r \leq b-1,$ we get
$p^{r+b+1}-\frac{p-1}{p^{b-r}}$ as an \textit{Average Order Limit}.
\begin{corollary}
Given any $p$-group $B$, and integer $r\geq2,$ there are sequences $G_n$, $H_n$ of $p$-groups with $|G_n|\to \infty$, $|H_n|\to \infty$ and $a(G_n \wr H_n)\to p^r a(B).$
\end{corollary}
\section{Proofs of the Theorems}
We start with two preliminary lemmas.
\begin{lemma}
For any semidirect product $H\rtimes K$, we have
$a(H\rtimes K)\geq a(K).$ Also, for $(h,k)\in H\rtimes K,$ we have order($k)|order((h,k)$).
\end{lemma}

\noindent {\bf Proof.} For any $(h,k)\in H\rtimes K,$ $(h,k)^m$ has
$k^m$ as 2nd coordinate. So, $(h,k)^m=1\implies k^m=1.$ So,
order($k)|order((h,k)$), in particular order($(h,k))\geq order(k$),
for all $h\in H.$ So, a($H\rtimes
K)\geq\frac{|H|}{|H|\cdot|K|}\sum_{k\in K}
order(k)$$=\frac{1}{|K|}\sum_{k\in K} order(k)=a(K).$

\begin{lemma}
Let $p$ be a prime and $b\geq 1$ an integer. Then
$a(\mathbf{Z}/p^b\mathbf{Z})=p^b+O(p^{b-1})$.
\end{lemma}

\noindent {\bf Proof.} In $\mathbf{Z}/p^b\mathbf{Z}$, there are
exactly $\phi(p^n)=p^n(1-p^{-1})$ elements of order $p^n$, for each
$1\leq n\leq b.$ So,
$$a(\mathbf{Z}/p^b\mathbf{Z}) =\frac{1+\sum_{n=1}^b
p^n\cdot p^n \cdot(1-p^{-1})}{p^b}=\frac{1+p(p-1)\sum_{n=0}^{b-1}
p^{2n}}{p^b}$$ $$=\frac{1+(p^2-p)\frac{p^{2b}-1}{p^2-1}}{p^b}=
\frac{1+\frac{p}{p+1}(p^{2b}-1)}{p^b}$$
$$=p^{-b}+\frac{p}{p+1}(p^b-p^{-b})=p^b+O(p^{b-1}).$$

Let us recall some standard facts before starting the proof of theorem 1.
 For groups $H$, $K$, and injective group homomorphism $\phi$: $K\rightarrow Aut(H)$, we have an injective group homomorphism
$\psi: H\rtimes_\phi K \rightarrow Perm(H)$ defined by
$\psi(h,k)(x)=h\phi(k)(x).$ So, $H\rtimes_\phi K$ can be regarded as
a subgroup of Perm($H$) consisting of the transformations $T_{h,k}$
defined by $T_{h,k} (x) = h\phi(k)(x).$

So, if $A$ and $B$ has at least $2$ elements, $A\wr B$ is the
subgroup of Perm($\prod_{b\in B} A$), consisting of the permutations
$T_{\underline{\alpha},x}$, $\underline{\alpha}\in\prod_{b\in B} A$,
$x\in B$, where
  $T_{\underline{\alpha},x}(\underline{a})=(\alpha_b a_{x^{-1}b})_b$, ($\underline{a}=(a_b)_b).$

  \noindent \vskip 3mm {\bf Proof of Theorem 1.}  Note that lemma 1 yields order($x)\bigm |$ order ($T_{\underline{\alpha}, x})$ for all $(\underline{\alpha}, x)\in A\wr B. $ Fix $x\in B.$ Let $d$=order($x$). For $m\geq1,$ we
   have $T_{\underline{\alpha}, x^{-1}}^m (\underline{a})=(\alpha_b \alpha_{xb} .. \alpha_{x^{m-1}b} a_{x^{m}b})_b.$

   So,
   $$T_{\underline{\alpha}, x^{-1}}^{dm}=1 \iff \alpha_b \alpha_{xb} .. \alpha_{x^{dm-1}b}=1~~\forall~~b\in
   B$$
$$ \iff (\alpha_b \alpha_{xb} .. \alpha_{x^{d-1}b})^m=1~~\forall~~b\in B.$$

Name this condition (1).

    Multiplication by $x$ divides $B$ into $\frac{|B|}{d}$ orbits of size $d$, and (1) is equivalent to saying that product of $\alpha_b$'s, for $b$ running in each orbit in cyclic order, has $m$'th power$=1$.

     Number of $\underline{\alpha}\in \prod_{b\in B} A$ satisfying (1) is exactly $(|A|^{d-1} s_m)^{\frac{|B|}{d}}.$ The reason is as follows. Fix a point $p$ in a orbit. For each $b\neq p$ in that orbit, $\alpha_b$ can be chosen to be anything in $A$. After that $\alpha_p$ has exactly $s_m$ choices. Altogether we get $\alpha_b$'s for $b$ running over a fixed orbit has exactly $|A|^{d-1} s_m$ choices. There are $\frac{|B|}{d}$ orbits. So in total we have $(|A|^{d-1} s_m)^{\frac{|B|}{d}}$ choices for $\underline{\alpha}$.

     From (1), it is clear that $T_{\underline{\alpha}, x^{-1}}^{d|A|}=1$. So, $\frac{order(T_{\underline{\alpha}, x^{-1}})}{d}  \bigm | |A|$

     For $m\geq1,$ let $g_m$= number of $\underline{\alpha}\in \prod_{b\in B} A$ with order($T_{\underline{\alpha}, x^{-1}})=dm$. We have just shown that unless $m\bigm| |A|,$ we have $g_m=0.$ Also, for all $k\geq1,$ $\sum_{m| k} g_m=$ number of $\underline{\alpha}$ with ($T_{\underline{\alpha},x^{-1}})^{dk}=1 = (|A|^{d-1} s_k)^{\frac{|B|}{d}}.$

     By Mobius inversion, $g_m=|A|^{\frac{|B|(d-1)}{d}} \sum_{n|m} \mu(\frac{m}{n})s_n^{\frac{|B|}{d}}$ for all $m\geq1.$ So, $\sum_{\underline{\alpha}\in \prod_{b\in B} A} order(T_{\underline{\alpha}, x^{-1}})$

     $=\sum_{m\bigm||A|} dmg_m=d|A|^{\frac{|B|(d-1)}{d}} \sum_{m\bigm||A|, n|m} m\mu(\frac{m}{n})s_n^{\frac{|B|}{d}}$

     $=d|A|^{\frac{|B|(d-1)}{d}} \sum_{n\bigm||A|, c\bigm| \frac{|A|}{n}} nc\mu(c)s_n^{\frac{|B|}{d}}$ ($m=nc$)

     $=d|A|^{\frac{|B|(d-1)}{d}} \sum_{m\bigm||A|, c\bigm| \frac{|A|}{m}} mc\mu(c)s_m^{\frac{|B|}{d}}$ (rename $n$ by $m$)

     $=d|A|^{\frac{|B|(d-1)}{d}} \sum_{m\bigm||A|} ms_m^{\frac{|B|}{d}}[\sum_{c\bigm| \frac{|A|}{m}} c\mu(c)]$

     $=d|A|^{\frac{|B|(d-1)}{d}} \sum_{m\bigm||A|} ms_m^{\frac{|B|}{d}}\tau(\frac{|A|}{m})$

     $=|A|^{|B|} d|A|^{-\frac{|B|}{d}} \sum_{m\bigm||A|} ms_m^{\frac{|B|}{d}}\tau(\frac{|A|}{m}).$

     This is true for any $x\in B$ or order $d$. There are $d_n$ elements of order $n$ in $B$, for each $n\bigm| |B|.$ So,

     $\sum_{x\in A\wr B} order(x)=\sum_{n\bigm| |B|} d_n |A|^{|B|} n|A|^{-\frac{|B|}{n}} [\sum_{m\bigm||A|} ms_m^{\frac{|B|}{n}}\tau(\frac{|A|}{m})] $

     $=|B| |A|^{|B|} \sum_{m\bigm||A|, n\bigm| |B|} d_n  \frac{n}{|B|} |A|^{-\frac{|B|}{n}} ms_m^{\frac{|B|}{n}}\tau(\frac{|A|}{m})$

     $=|B| |A|^{|B|} \sum_{m\bigm||A|, n\bigm| |B|} d_{\frac{|B|}{n}}  n^{-1} |A|^{-n} ms_m^n \tau(\frac{|A|}{m})$ ( by replace $n$ with $\frac{|B|}{n}$)

     $=|B| |A|^{|B|} \sum_{m\bigm| |A|, n\bigm| |B|} \frac{m}{n} (\frac{s_m}{|A|})^n d_\frac{|B|}{n} \tau(\frac{|A|}{m})$.

 Dividing by $|B| |A|^{|B|}$ we get the desired result.

  \noindent \vskip 3mm {\bf Proof of Theorem 2.}
  Theorem 1, together with the observation that $s_{p^m}=p^a$ for $m\geq d$ yields

  a($A\wr B$)

  = $\sum_{n\leq b} \frac{p^a}{p^n} (\frac{s_{p^a}}{p^a})^{p^n} d_{p^{b-n}}$

  $-(p-1)\sum_{m\leq a-1, n\leq b} \frac{p^m}{p^n}(\frac{s_{p^m}}{p^a})^{p^n} d_{p^{b-n}}$

  $= p^{a-b}\sum_{n\leq b} p^{b-n} d_{p^{b-n}}$

  $-(p-1)\sum_{d\leq m \leq a-1, n\leq b} \frac{p^m}{p^n} d_{p^{b-n}}$

  $-(p-1)\sum_{m\leq d-1, n\leq b} \frac{p^m}{p^n}(\frac{s_{p^m}}{p^a})^{p^n} d_{p^{b-n}}$

  $=p^a$a($B$)$-(p-1)p^d[\sum_{m\leq a-d-1} p^m] \sum_{n\leq b} p^{-n} d_{p^{b-n}}$

  $-(p-1)\sum_{m\leq d-1, n\leq b} \frac{p^m}{p^n}(\frac{s_{p^m}}{p^a})^{p^n} d_{p^{b-n}}$

  $=p^a$a($B$)$-(p^{a-d}-1)p^d .a(B)$

  $-(p-1)\sum_{m\leq d-1, n\leq b} \frac{p^m}{p^n}(\frac{s_{p^m}}{p^a})^{p^n} d_{p^{b-n}}$

  $=p^d$a($B$) $-(p-1)\sum_{m\leq d-1, n\leq b} p^m(\frac{s_{p^m}}{p^a})^{p^n} p^{-n} d_{p^{b-n}}$

$=p^d$a($B$) $-(p-1)\sum_{n\leq b} p^{-n} d_{p^{b-n}} [\sum_{m\leq d-1} p^m(\frac{s_{p^m}}{p^a})^{p^n}] $

$=p^d$a($B$) $-(p-1)a(B)\sum_{n\leq b} k_n. [\sum_{m\leq d-1} p^m(\frac{s_{p^m}}{p^a})^{p^n}] $.

 \noindent \vskip 3mm {\bf Proof of theorem 3.} Theorem 2 proves the second inequality, and the first inequality follows from lemma 1.

  \noindent \vskip 3mm {\bf Proof of theorem 4.} Let us look at theorem 2 more closely. $\frac{s_{p^m}}{p^a}<1$ $\forall m\leq d-1$. So, $\sum_{m\leq d-2} p^m (\frac{s_{p^m}}{p^a})^{p^n}\leq \sum_{m\leq d-2} p^m=\frac{p^{d-1}-1}{p-1},$ for each $n\leq b$. We also have $k_n\geq 0$ $\forall n$, $\sum_{n\leq b} k_n=1.$ Hence,

  $(p-1)a(B).\sum_{n\leq b}k_n . [\sum_{m\leq d-2} p^m (\frac{s_{p^m}}{p^a})^{p^n}]=O(p^{d-1} a(B)).$

  \noindent \vskip 3mm {\bf Proof of theorem 5.}
    Let $x\in B$, order ($x$)=$p^{e_1}$, $\underline{\alpha}\in \prod_{b\in B} A.$
Multiplication by $x$ divides $B$ into $p^{b-e_1}$ orbits of size
$p^{e_1};$ let $b_1, b_2,..., b_{p^{b-e_1}}$ be representatives of
distinct orbits. By (1) of theorem 1, the order of
$T_{\underline{\alpha},x^{-1}}$ equals $p^{e_1}.\max_{1\leq i \leq
p^{b-e_1}} order(\alpha_{b_i}\alpha_{xb_i} \cdots
\alpha_{x^{p^{e_1}-1}b_i}) \leq p^e.p^d=p^{d+e}.$

    If $x\in B$ is of order $p^e,$ and $y\in A$ is of order $p^d$, then
the order of $T_{\underline{\alpha},x^{-1}}$ equals $p^{d+e},$
where $\underline{\alpha}$ is defined by $\alpha_1=y,$ $\alpha_b=1$
$\forall b\neq 1.$ So, m($A \wr B)=p^{d+e}.$

    Now we prove the second statement.

    Again fix $x\in B$, order ($x$)=$p^{e_1}$. Let $k\leq d+e-e_1.$ Note that
    $$T_{\underline{\alpha},x^{-1}}^{p^{d+e-k}}=1 \iff
(\alpha_{b_i}\alpha_{xb_i}...\alpha_{x^{p^{e_1}-1}b_i})^{p^{d+e-e_1-k}}=1~~\forall~~i.$$

Name this condition (2).

    As we showed in the proof of theorem 1, the number of $\underline{\alpha}\in \prod_{b\in B} A$ satisfying (2)
is
$((p^a)^{p^{e_1}-1}s_{p^d+e-e_1-k})^{p^{d-e_1}}=(p^{ap^{e_1}}r_{A,k-(e-e_1)})^{p^{b-e_1}}=p^{ap^b}r_{A,k-(e-e_1)}^{p^{b-e_1}}$.
(Note that $s_{p^d+e-e_1-k}=p^a r_{A,k-(e-e_1)}.$)

    So, for $k\leq d+e-e_1,$ $|\{\underline{\alpha}\in \prod_{b\in B} A: T_{\underline{\alpha},x^{-1}}^{p^{d+e-k}}=1\}|=p^{ap^b}r_{A,k-(e-e_1)}^{p^{b-e_1}}. $ This is true for $k> d+e-e_1$ also, as then both sides are $0$.

    So, $r_{A\wr B,k}$

    $=p^{-ap^b}p^{-b}\sum_{e_1=0}^e(\text{no. of elements of $B$ of order $p^{e_1}$})p^{ap^b}r_{A,k-(e-e_1)}^{p^{b-e_1}}$

    $=\sum_{i=0}^e \frac{\text{no. of elements of $B$ of order $p^{e-i}$}}{p^b}\cdot r_{A,k-i}^{p^{b-e+i}}$ (here $i=e-e_1$)

    $=\sum_{i=0}^e (r_{B,i}-r_{B,i+1}) r_{A,k-i}^{p^{b-e+i}}$.

    Corollary 5.1 is a direct consequence of theorem 5.

    To prove Corollary 5.2, we need a lemma. Write $r_{n,k}=r_{A_n,k}.$

    \begin{lemma}
    For all $k$, $\lim_{n\to \infty} r_{n,k}=1.$
    \end{lemma}

    \noindent {\bf Proof.} We proceed by induction on $k$. For $k\leq 0,$ we have $r_{n,k}=1$ $\forall n,$ and nothing to show. For the induction step, let $k\geq 1$ and assume $\lim_{n\to \infty} r_{n,k-1}=1.$ By corollary 5.1., $r_{n+1,k}=(1-p^{-1})r_{n,k}+\frac{r_{n,k-1}^p}{p}.$ So,

    $\liminf_{n\to \infty} r_{n+1,k}=(1-p^{-1})\liminf_{n\to \infty} r_{n,k}+p^{-1}$,
    as induction hypothesis implies $\lim_{n\to \infty}
    r_{n,k-1}=1.$

    But $\liminf_{n\to \infty} r_{n+1,k}=
\liminf_{n\to \infty} r_{n,k}$. Hence, $\liminf_{n\to \infty}
r_{n,k}=(1-p^{-1})\liminf_{n\to \infty} r_{n,k}+ p^{-1}$; that is,
$\liminf_{n\to \infty} r_{n,k}=1$. Since we always have $0\leq
r_{n,k}\leq 1,$ we get $\lim_{n\to \infty} r_{n,k}=1$. Induction
completes the proof.

\noindent \vskip 3mm {\bf Proof of corollary 5.2.} Let $m(A)=p^d$,
$|B|=p^b$. By theorem 5, $m(A_n)=p^{d+n}$ for all $n$. By theorem 2,

$\psi(A_n,B)= 1-\frac{p-1}{p^{d+n}}\sum_{r\leq b} k_r (\sum_{m=0}^{d+n-1} p^m r_{n,d+n-m}^{p^r})$

$=1-(1-p^{-1})\sum_{r\leq b} k_r (\sum_{m=0}^{d+n-1} p^{-(d+n-1-m)} r_{n,(d+n-1-m)+1}^{p^r})$

$= 1-(1-p^{-1})\sum_{r\leq b} k_r (\sum_{m=0}^{d+n-1} p^{-m} r_{n,m+1}^{p^r})$ (replace $m$ by $d+n-1-m$).

Fix $M\in \mathbf{N}.$ For all sufficiently large $n$, we have $d+n-1>M$, hence $\psi(A_n,B)\leq 1-(1-p^{-1})\sum_{r\leq b} k_r (\sum_{m=0}^{M} p^{-m} r_{n,m+1}^{p^r}).$ Since $\lim_{n\to \infty} r_{n,m+1}=1$ by lemma 3, we get

$\limsup_{n\to \infty} \psi(A_n,B)\leq 1-(1-p^{-1})\sum_{r\leq b} k_r \sum_{m=0}^{M} p^{-m}$

$=1-(1-p^{-1}) \sum_{m=0}^{M} p^{-m}$ (as $\sum_{r\leq b} k_r=1$) for all $M\in \mathbf{N}.$

Taking $M\to \infty$, we get $\limsup_{n\to \infty} \psi(A_n,B) \leq 0$, that is, $\lim_{n\to \infty} \psi(A_n,B) = 0.$

\noindent \vskip 3mm {\bf Proof of theorem 6.} We use the notation
of theorem 2, and write $t$ instead of $t(A)$. Projection gives a
surjective group homomorphism $\phi: A\to
(\mathbf{Z}/p^d\mathbf{Z})^{t}$. Note that the subgroup of
$p^{d-1}$-torsion elements of $A$ is
$\phi^{-1}((p\mathbf{Z}/p^d\mathbf{Z})^{t})$, and for $0\leq m\leq
d-2,$ the subgroup of $p^{d-m}$-torsion elements of $A$ is contained
in $\phi^{-1}((p^m \mathbf{Z}/p^d\mathbf{Z})^{t})$. So,
$\frac{s_{p^{d-1}}}{p^a}=p^{-t}$, and for each $0\leq m\leq d-2,$ we
have $\frac{s_{p^m}}{p^a}\leq p^{-mt}\leq p^{-2t}.$ By theorem 2,
$$p^d a(B)-a(A\wr B)-(p-1)a(B)k_0\cdot [\sum_{m\leq d-1}
p^m\frac{s_{p^m}}{p^a}]$$
$$=(p-1)a(B)\sum_{1\leq n\leq b} k_n\cdot
[\sum_{m\leq d-1} p^m(\frac{s_{p^m}}{p^a})^{p^n}]$$ $$\leq
(p-1)a(B)\sum_{1\leq n\leq b} k_n \sum_{m\leq d-1}
\frac{p^m}{p^{tp^n}}$$ $$=(p-1)a(B)(\sum_{1\leq n\leq b} k_n
p^{-tp^n})\cdot\frac{p^d-1}{p-1}$$ $$\leq a(B)\sum_{1\leq n\leq b}
p^{d-tp^n}k_n$$ $$\leq a(B)p^{d-tp}\sum_{1\leq n \leq b} k_n \leq
a(B)p^{d-tp}.$$

If $B$ is noncyclic, $k_0=0;$ so $p^d a(B)-a(A\wr B) \leq
a(B)p^{d-tp}.$ Dividing by $p^d a(B)$, we get $1-p^{-tp}\leq
\psi(A,B).$ If $B$ is the cyclic group $\mathbf{Z}/p^b\mathbf{Z}$,
we have $k_0=\frac{\phi(p^b)}{a(B)}=\frac{p-1}{p}\cdot
\frac{p^b}{a(B)}=(1+O(\frac{1}{p}))\cdot(1+O(\frac{1}{p}))=1+O(\frac{1}{p}).$
Here, we used lemma 2, and the observation that
$(1+O(p^{-1}))^{-1}=1+O(p^{-1}).$ So, $\sum_{m\leq d-1}
p^m\frac{s_{p^m}}{p^a}=p^{d-1-t}+O(\sum_{m\leq d-2} p^{m-2t})=
p^{d-1-t} + O(\frac{p^{d-1}-1}{p-1}\cdot p^{-2t})= p^{d-1-t} +
O(p^{d-2-2t})$.

Hence, $\frac{(p-1)a(B)k_0\cdot [\sum_{m\leq d-1}
p^m\frac{s_{p^m}}{p^a}]}{p^d
a(B)}=\frac{(p-1)(1+O(p^{-1}))(p^{d-1-t}+O(p^{d-2-2t}))}{p^d}=\frac{p^{d-t}+O(p^{d-1-t})}{p^d}=p^{-t}+O(p^{-t-1}).$

 We have shown

 $0\leq p^d a(B)-a(A\wr B)-(p-1)a(B)k_0\cdot [\sum_{m\leq d-1} p^m\frac{s_{p^m}}{p^a}]\leq a(B)p^{d-tp}.$

 So, dividing by $p^d a(B)$ we get

 $\psi(A,B)=1-\frac{(p-1)a(B)k_0\cdot [\sum_{m\leq d-1} p^m\frac{s_{p^m}}{p^a}]}{p^d a(B)}$ $+O(p^{-tp})$

 $=1+p^{-t}+O(p^{-t-1})+O(p^{-tp})=1+p^{-t}+O(p^{-t-1}).$

\noindent \vskip 3mm {\bf Proof of theorem 7.} Let $|B_n|=p^{b_n}$.
By theorem 2,
$$a((\mathbf{Z}/p\mathbf{Z})^n \wr
B_n) = p \cdot a(B_n)-(p-1)a(B_n) \sum_{m\leq b_n}
k_{m,n}(\frac{1}{p^n})^{p^m},$$ for some $0\leq k_{m,n}\leq 1,$ with
$\sum_{m\leq b_n} k_{m,n}=1$ for each $n.$

So, for all sufficiently large $n$, so that $a(B_n)\leq \beta+1;$ we
have $$|a((\mathbf{Z}/p\mathbf{Z})^n\wr B_n)-p\beta|\leq
p|a(B_n)-\beta|+(p-1)(\beta+1)p^{-n} \to 0.$$

\noindent \vskip 3mm {\bf Proof of corollary 7.1.} Let
$G_n=(\mathbf{Z}/p\mathbf{Z})^n,$
$H_n=(\mathbf{Z}/p\mathbf{Z})^n\wr((\mathbf{Z}/p\mathbf{Z})^n\wr(...((\mathbf{Z}/p\mathbf{Z})^n\wr
B))),$ there are $r-1$  $(\mathbf{Z}/p\mathbf{Z})^n$'s here. Now
corollary 7.1 follows by repeated application of theorem 7.

\section{Acknowledgement}
I am grateful to Professor B. Sury of Indian Statistical Institute,
Bangalore, who encouraged me to investigate many of the questions I
addressed in this article. He also helped me very much to arrange
the article in proper order. \vskip 5mm
\section{Declaration of interest}
No potential conflict of interest was reported by the author.

\end{document}